\documentclass{llncs}
\usepackage[utf8]{inputenc}
\usepackage[english]{babel}
\usepackage{amsmath}
\usepackage{amsfonts}
\usepackage{amssymb}
\usepackage{graphicx}
\usepackage{color}
\usepackage{textcomp}
\usepackage{comment}

\newcommand{\R}{\mathbb{R}}
\newcommand{\Hb}{\mathbb{H}}
\newcommand{\Sp}{\mathbb{S}}
\newcommand{\radon}[2]{R_{#1}^{}}
\newcommand{\eqdef}{\overset{def}{=}}	

\title{An analog of Chang inversion formula for weighted Radon transforms in multidimensions\thanks{The main part of the work was fulfilled during the stage of the first author in the Centre 
de Math\'ematiques Appliqu\'ees of Ecole Polytechnique in March-May 2016.}}
\author{Fedor O. Goncharov\inst{1}\textsuperscript{,} \inst{2}, 
Roman G. Novikov\inst{3}\textsuperscript{,}\inst{4}}
\institute{Moscow Institute of Physics and Technology, \\9 Institutskiy per., 
Dolgoprudny, Moscow Region, 141700, Russian Federation;\\
\and Institute for Information Transmission Problems (RAS),\\ 
Bolshoy Karetny per. 19, buld.1, Moscow 127051, Russian Federation;\\
\and CNRS (UMR 7641), Centre de Math\'ematiques Applique\'es, Ecole Polytechnique,\\
91128 Palaiseau, France;\\
\and Institute of Earthquake Prediction Theory and Mathematical Geophysics (RAS),\\ 
117997 Moscow, Russian Federation;\\
\email{fedor.goncharov.ol@gmail.com, novikov@cmap.polytechnique.fr}}

\begin{document}

\maketitle
\begin{abstract}
In this work we study weighted Radon transforms in multidimensions. We introduce 
an analog of Chang approximate inversion formula for such transforms and describe all
weights for which this formula is exact. In addition, we indicate possible tomographical
applications of inversion methods for weighted Radon transforms in 3D. 
\end{abstract}

\section{\large Introduction}
We consider the weighted Radon transforms $\radon{W}{d,d-1}$ defined by the formula
\begin{align}\label{radon_general}
	&\radon{W}{d,d-1}f(s,\theta) \eqdef \int\limits_{x\theta = s} W(x,\theta) f(x) dx,	\\
	&(s,\theta) \in \R\times \Sp^{n-1}, \, x\in \R^n,\, n \geq 2, \nonumber
\end{align}
where $W = W(x,\theta)$ is the weight, $f = f(x)$ is a test function; see e.g. \cite{boman1987support}.
Such transforms arise in many domains of pure and applied mathematics; see e.g. \cite{deans2007radon}, \cite{durrani1984radon},
\cite{gel2014integral}, \cite{grangeat1991mathematical}, \cite{kunyansky1992generalized}, \cite{natterer1986mathematics}.
In the present work we assume that 
\begin{align}\label{weight_class}
	&W \text{ is complex -- valued},\nonumber\\
	&W \in C(\R^n \times \Sp^{n-1}) \cap L^{\infty}(\R^n \times \Sp^{n-1}),\\
	&w_0(x) \eqdef  \dfrac{1}{|\Sp^{n-1}|}\int\limits_{\Sp^{n-1}}W(x,\theta)d\theta \neq 0,\, x\in \R^n, \nonumber
\end{align}
where $d\theta$ is the element of standard measure on $\Sp^{n-1}$, $|\Sp^{n-1}|$ is the standard measure of $\Sp^{n-1}$.

\par If $W \equiv 1$, then $R= \radon{W}{d,d-1}$ is the classical Radon transform in $\R^n$; see for example  \cite{gel2014integral}, \cite{helgason2011radon}, \cite{ludwig1966radon}, \cite{radon1994uber}. Explicit inversion formulas for
$R$ were given for the first time in \cite{radon1994uber}.\\
\indent In dimension $n = 2$, the transforms $\radon{W}{d,d-1}$ are also known as weighted ray transforms on the plane; see e.g. \cite{kunyansky1992generalized}, \cite{natterer1986mathematics}. For several important cases of $W$ satisfying \eqref{weight_class} for $d=2$, explicit (and exact) inversion formulas for $\radon{W}{}$ were obtained in \cite{boman2004novikov}, \cite{gindikin2010remark}, \cite{novikov2002inversion}, \cite{novikov2011weighted}, \cite{tretiak1980exponential}.\\
\indent On the other hand, it seems that no explicit inversion formulas for $\radon{W}{}$ were given yet in the literature under assumptions 
\eqref{weight_class} for $n\geq 3$, if $W\neq w_0$.\\
\indent In the present work we introduce an analog of Chang approximate (but explicit) inversion formula for $R_W$ under assumptions \eqref{weight_class}, for $n\geq 3$, and describe all $W$ for which this formula is exact. These results are presented in Section 2.\\
\indent In addition, we indicate possible tomographical applications of inversion methods for 
$R_W$ in dimension $n=3$. These considerations are presented in Section 3.

\renewcommand{\thesection}{\large 2}
\section{\large Chang-type formulas in multidimensions }
We consider the following approximate inversion formulas for $R_W$ under assumptions \eqref{weight_class} in dimension $n\geq 2$:
\begin{align}\label{chang}
	&f_{appr}(x) \eqdef \dfrac{(-1)^{(n-2)/2}}{2(2\pi)^{n-1} w_0(x)}\int\limits_{\Sp^{n-1}} \mathbb{H} \left[ R_Wf \right]^{(n-1)}(x\theta,\theta)d\theta, \\
	& x\in \R^n, \, n\text{ is even}, \nonumber\\
	\label{chang_2}
	&f_{appr}(x) \eqdef \dfrac{(-1)^{(n-1)/2}}{2(2\pi)^{n-1} w_0(x)}
	\int\limits_{\Sp^{n-1}} \left[R_Wf\right]^{(n-1)}(x\theta, \theta)d\theta,\\ 
	& x\in \R^n, \, n\text{ is odd},\nonumber\\
	&\left[R_Wf\right]^{(n-1)}(s,\theta) = \dfrac{d^{n-1}}{ds^{n-1}} \radon{W}{d,d-1}f(s,\theta), s\in \R, \, \theta\in \Sp^{n-1},\\
	&\mathbb{H}\phi(s) \eqdef \dfrac{1}{\pi} p.v. \int\limits_{\R} \dfrac{\phi(t)}{s-t}dt, \, s\in \R.
\end{align}
\indent For $W\equiv 1$ formulas \eqref{chang}, \eqref{chang_2} are exact, i.e. $f_{appr} = f$, and are known as the classical Radon inversion formulas, going back to \cite{radon1994uber}. 
\par As a corollary of the classical Radon inversion formulas and definition \eqref{radon_general}, formulas \eqref{chang}, \eqref{chang_2} for $W\equiv w_0$ are also exact.\\
\indent Formula \eqref{chang} for $n=2$ is known as Chang approximate inversion formula for weighted Radon transforms on the plane. This explicit but approximate inversion formula was suggested for the first time in \cite{chang1978method} for the case when
\begin{align}\label{weight_1}
		&W(x,\theta) = \exp\left( -Da(x,\theta^\perp)\right),\\
		\label{weight_2}
		&Da(x,\theta^\perp) = \int\limits_{0}^{+\infty} a(x + t\theta^\perp) dt,
\end{align}
where $a$ is a non-negative sufficiently regular function on $\R^2$ with compact support, and $\theta = (\theta_1,\theta_2)\in \Sp^{n-1}, \, \theta^\perp = (\theta_2, -\theta_1)$. We recall that $R_W$ for $W$ given by \eqref{weight_1}, \eqref{weight_2} is known as attenuated Radon transform on the plane and arises, in particular, in the single photon emission tomography (SPECT). In this case an explicit and simultaneously exact inversion formula for $R_W$ was obtained for the first time in \cite{novikov2002inversion}.\\
\indent We emphasize that formulas \eqref{chang}, \eqref{chang_2} are approximate, in general. In addition, the following result holds:
\begin{theorem}\label{main_theorem}
	Let $W$ satisfy \eqref{weight_class}. Let $f_{appr}$ be defined by \eqref{chang}, \eqref{chang_2} in terms of $R_Wf$ and $w_0$, $n\geq 2$.
	Then $f_{appr}=f$  (in the sense of distributions) on $\R^n$ for all $f\in C_0(\R^n)$ if and only if 
	\begin{equation}\label{condition}
		W(x,\theta) - w_0(x) \equiv w_0(x) - W(x,-\theta), \, x\in \R^n, \, \theta\in \Sp^{n-1}.
	\end{equation}
\end{theorem}
Here $C_0(\R^n)$ denotes the space of all continous compactly supported functions on $\R^n$.\\
\indent The result of Theorem 1 for $n=2$ was obtained for the first time in \cite{novikov2011weighted}. Theorem 1 in the general case is proved in Section 4.

If $W$ satisfy \eqref{weight_class}, $f\in C_0(\R^n)$, but the the symmetry condition \eqref{condition} does not hold, i.e.
	$$
		w_0(x) \neq \dfrac{1}{2}\left( W(x,\theta) + W(x,-\theta) \right), \text{ for some } x\in \R^n, \, \theta\in \Sp^{n-1},
	$$
	then \eqref{chang}, \eqref{chang_2} can be considered as approximate formulas for finding $f$ from $R_Wf$.

\renewcommand{\thesection}{\large 3}
\section{\large Weighted Radon transforms in 3D in tomographies}
In several tomographies the measured data are modeled by weighted ray transforms $P_wf$ defined by the formula
\begin{align}\label{weighted_ray_transform}
	&P_wf(x, \alpha) = 
	\int\limits_{\R} w(x + \alpha t, \alpha) f(x + \alpha t)\, dt, \, 
	(x, \alpha)\in T\Sp^2, \\
	& T\Sp^{2} = \{(x,\alpha) \in \R^3 \times \Sp^2 : x\alpha = 0\}, \nonumber
\end{align}
where $f$ is an object function defined on $\R^3$, $w$ is the weight function defined on $\R^3\times \Sp^2$, and $T\Sp^2$ can be considered as the set of all rays (oriented straight lines) in $\R^3$. In particular, in the case of the single-photon emission computed tomography (SPECT) 
the weight $w$ is given by formulas \eqref{weight_1}, \eqref{weight_2}, where $\theta^\perp = \alpha \in 
\Sp^{2}, \, x\in \R^3$.
\par In practical tomographical considerations $P_wf( x,\alpha)$ usually arises for rays $(x,\alpha)$ parallel to some fixed plane 
\begin{equation}
	\Sigma_{\eta} = \{x\in \R^3 : x\eta = 0\}, \, \eta\in \Sp^{2},
\end{equation}
i.e., for $\alpha\eta = 0$.

\par The point is that the following formulas hold:
\begin{align}\label{averaging_1}
	&R_Wf(s,\theta) = \int\limits_{\R}P_wf(s\theta + \tau [\theta, \alpha], \alpha)d\tau, \, s\in \R,
	\, \theta\in\Sp^2,\\
	\nonumber
	&W(x,\theta) = w(x,\alpha),\, \, \alpha = \alpha(\eta, \theta) = \dfrac{[\eta, \theta]}{|[\eta, \theta]|},
	 \, \, [\eta,\theta]\neq 0, \, \, x\in \R^3,
\end{align}
where $[\cdot, \cdot]$ stands for the standart vector product in $\R^3$.
\par Due to formula \eqref{averaging_1} the measured tomographical data modeled by $P_wf$ can
be reduced to averaged data modeled by $R_Wf$. In particular, this reduction drastically reduces
the level of random noise in the initial data.
\par Therefore, formula \eqref{chang_2} for $n=3$ and other possible methods for finding $f$ from
$R_Wf$ in 3D may be important for tomographies, where measured data are modeled by $P_wf$ of 
\eqref{weighted_ray_transform}.
\begin{remark}
	The weight $W$ arising in \eqref{averaging_1} is not continuous, in general. However, the result of Theorem~\ref{main_theorem} remains valid for this $W$, at least, under the assumptions that
	$w$ is bounded and continuous on $\R^3\times \Sp^2$, and $w_0(x)\neq 0, \, x\in \R^3$, where 
	$w_0$ is defined in \eqref{weight_class}.
\end{remark}

\renewcommand{\thesection}{\large 4}
\section{\large Proof of Theorem 1}
For $W$ satisfying \eqref{weight_class} we also consider its symmetrization defined by
\begin{equation}\label{symm_weight}
	W_s(x,\theta) \eqdef \dfrac{1}{2}\left(W(x,\theta) + W(x,-\theta)\right), \, 
	x\in \R^n, \, \theta\in \Sp^{n-1}.
\end{equation}
Using definitions \eqref{radon_general}, \eqref{symm_weight} we obtain 
\begin{equation}\label{radon_symmetry}
	R_{W_{s}}f(s,\theta) = \dfrac{1}{2}
	\left(
		R_Wf(s,\theta) + R_Wf(-s,-\theta)
	\right).
\end{equation}
In addition, if $W$ satisfies \eqref{condition}, then 
\begin{equation}\label{weight_symmetry}
	W_s(x,\theta) = w_0(x), \, x\in \R^n, \, \theta\in \Sp^{n-1}.
\end{equation}

\subsection{Proof of sufficiency}
The sufficiency of symmetry \eqref{condition} follows from formulas \eqref{chang}, \eqref{chang_2} for the exact case with $W\equiv w_0$, the identities
\begin{align}\label{even_identity}
	&f_{appr}(x) = \dfrac{(-1)^{(n-2)/2}}{2(2\pi)^{n-1} w_0(x)}\int\limits_{\Sp^{n-1}} \mathbb{H} \left[ R_{W_s}f \right]^{(n-1)}(x\theta,\theta)d\theta, \\
	&\text{for even }n,\nonumber\\ 
	\label{odd_identity}
	&f_{appr}(x) = \dfrac{(-1)^{(n-1)/2}}{2(2\pi)^{n-1} w_0(x)}\
	\int\limits_{\Sp^{n-1}} \left[R_{W_s}f\right]^{(n-1)}(x\theta, \theta)d\theta,\\
	&\text{for odd }n, \nonumber
\end{align}
 and from the identities \eqref{radon_symmetry}, \eqref{weight_symmetry}.
\par In turn, \eqref{even_identity} follows from the identities
\begin{align}
	&\int\limits_{\Sp^{n-1}} \mathbb{H}\left[
		R_{W}f
	\right]^{(n-1)}(x\theta,\theta)d\theta\nonumber\\ \label{hilbert_symmetry}
	&= \dfrac{1}{2}
		\int\limits_{\Sp^{n-1}} \left(\mathbb{H}\left[
			R_{W}f
		\right]^{(n-1)}(x\theta,\theta) + 
								   \mathbb{H}
								   \left[
			R_{W}f
		\right]^{(n-1)}(-x\theta,-\theta)
		\right)
	d\theta \\
	&=\int\limits_{\Sp^{n-1}}\mathbb{H}\left[R_{W_s}f\right]^{(n-1)}(x\theta,\theta)d\theta.\nonumber
\end{align}
In addition, the second of the identities of \eqref{hilbert_symmetry} follows from the identities:
\begin{align}
	&\mathbb{H}\left[
		R_{W_s}f
	\right]^{(n-1)}(s,\theta) = \dfrac{1}{2\pi}p.v.\int\limits_{\R}\dfrac{1}{s-t}\times\nonumber\\
	&\qquad \times \dfrac{d^{n-1}}{dt^{n-1}}
	\Big[R_Wf(t,\theta) + 
	R_Wf(-t,-\theta)\Big]dt\nonumber\\
	&\qquad = \dfrac{1}{2}\mathbb{H}\Big[ R_Wf\Big]^{(n-1)}(s,\theta) + \dfrac{(-1)^{n-1}}{2\pi}p.v.
	\int\limits_{\R}\dfrac{\left[R_Wf\right]^{(n-1)}(-t,-\theta)}{s-t}dt;\\
	&\dfrac{(-1)^{n-1}}{\pi}p.v.
	\int\limits_{\R}\dfrac{\left[R_Wf\right]^{(n-1)}(-t,-\theta)}{s-t}dt = - \dfrac{(-1)^{n-1}}{\pi}p.v.
	\int\limits_{\R}\dfrac{\left[R_Wf\right]^{(n-1)}(t,-\theta)}{-s-t}dt \nonumber\\
	&\qquad =(-1)^n\mathbb{H}\left[R_Wf \right]^{(n-1)}(-s,-\theta) = 
	\mathbb{H}\left[R_Wf \right]^{(n-1)}(-s,-\theta).
\end{align}
This concludes the proof of sufficiency for $n$ even.\\
\indent Finally, \eqref{odd_identity} follows from the identities
\begin{align}
	\int\limits_{\Sp^{n-1}}\left[R_Wf\right]^{(n-1)}&(x\theta,\theta)d\theta \nonumber\\
	\label{sufficiency_odd_1}
	=\dfrac{1}{2}\int\limits_{\Sp^{n-1}}&\left(\left[R_Wf\right]^{(n-1)}(x\theta,\theta) + 
	\left[R_Wf\right]^{(n-1)}(-x\theta,-\theta)\right)d\theta,\\ \label{sufficiency_odd_2}
	\left[R_{W_s}f\right]^{(n-1)}(t,\theta) &= \dfrac{1}{2}\dfrac{d^{n-1}}{dt^{n-1}}
	\Big[
		\left[R_Wf\right](t,\theta) + \left[R_Wf\right](-t,-\theta)
	\Big] \nonumber\\
	&= \dfrac{1}{2}\Big[
		\left[R_Wf\right]^{(n-1)}(t,\theta) + (-1)^{n-1}\left[R_Wf\right]^{(n-1)}(-t,-\theta)
	\Big]\nonumber\\
	&= \dfrac{1}{2}\Big[
		\left[R_Wf\right]^{(n-1)}(t,\theta) + \left[R_Wf\right]^{(n-1)}(-t,-\theta)
		\Big]. 
\end{align}
This concludes the proof of sufficiency for odd $n$.

\subsection{Proof of necessity}
Using that $f_{appr} = f$ for all $f\in C_0(\R^n)$ and using formulas \eqref{chang}, \eqref{chang_2} for the exact case $W\equiv w_0$, we obtain
\begin{align}
	\label{condition_necessity_even}
	&\int\limits_{\Sp^{n-1}} \left(\mathbb{H}\left[R_{W}f\right]^{(n-1)}(x\theta,\theta) - 
	\mathbb{H}\left[R_{w_0}f\right]^{(n-1)}(x\theta,\theta)\right)d\theta = 0\\
	&\text{on }\R^n\text{ for even }n,\nonumber\\
	\label{condition_necessity_odd}
	&\int\limits_{\Sp^{n-1}}	\left[R_{W}f - R_{w_0}f\right]^{(n-1)}
	(x\theta,\theta)d\theta = 0\\
   &\text{on }\R^n\text{ for odd }n, \nonumber
\end{align}
for all $f\in C_0(\R^n)$.
\par Identities \eqref{hilbert_symmetry}, \eqref{sufficiency_odd_1}, \eqref{sufficiency_odd_2}, \eqref{condition_necessity_even}, \eqref{condition_necessity_odd} imply the identities 
\begin{align}
    \label{necessity_even_sym}
	&\int\limits_{\Sp^{n-1}} \left(\mathbb{H}\left[R_{W_s}f\right]^{(n-1)}(x\theta,\theta) - 
	\mathbb{H}\left[R_{w_0}f\right]^{(n-1)}(x\theta,\theta)\right)d\theta = 0\\
	&\text{on } \R^n\text{ for even }n,\nonumber\\
	\label{necessity_odd_sym}
	&\int\limits_{\Sp^{n-1}}	\left[R_{W_s}f - R_{w_0}f\right]^{(n-1)}
	(x\theta,\theta)d\theta = 0\\
   &\text{on }\R^n \text{ for odd }n, \nonumber
\end{align}
for all $f\in C_0(\R^n)$.
\par The necessity of symmetry \eqref{condition} follows from the identities \eqref{necessity_even_sym}, \eqref{necessity_odd_sym} and the following lemmas:
\begin{lemma}
	Let \eqref{necessity_even_sym}, \eqref{necessity_odd_sym} be valid for fixed $f\in C_0(\R^n)$ and $W$ satisfying \eqref{weight_class}, $n\geq 2$. Then 
	\begin{equation}
		R_{W_s}f = R_{w_0}f. \label{radon_eq}
	\end{equation}
\end{lemma}
\begin{lemma}
	Let \eqref{radon_eq} be valid for all $f\in C_0(\R^n)$ and fixed $W$ satisfying \eqref{weight_class}, $n\geq 2$. Then
	\begin{equation}
		W_s = w_0.
	\end{equation}
\end{lemma}
Lemmas 1 and 2 are proved in Sections 5 and 6.

\renewcommand{\thesection}{\large 5}
\section{\large Proof of Lemma 1}
We will use the following formulas
\begin{align}\label{fourier_dual_radon}
	&\int\limits_{\R^n}e^{i\xi x}\nonumber
	\int\limits_{\Sp^{n-1}}g(x\theta,\theta)\, d\theta \, dx  \\
	&=\dfrac{\sqrt{2\pi}}{|\xi|^{n-1}}\left(
		\hat{g}\left(
			|\xi|, \dfrac{\xi}{|\xi|}
			\right)	+
		\hat{g}\left(
			-|\xi|, -\dfrac{\xi}{|\xi|}
			\right)			
	\right), \\ \label{fourier_dual_derivative}
	&\int\limits_{\R^n}e^{i\xi x}\nonumber
	\int\limits_{\Sp^{n-1}}g^{(n-1)}(x\theta,\theta)\, d\theta\, dx =
	\int\limits_{\R^n}e^{i\xi x}\nonumber
	\int\limits_{\Sp^{n-1}}(\theta\nabla_x)^{n-1}g(x\theta,\theta)\, d\theta \, dx\\
	&=(-i)^{n-1}\sqrt{2\pi}
	\left(
		\hat{g}\left(
			|\xi|, \dfrac{\xi}{|\xi|}
			\right)	+ (-1)^{n-1}
		\hat{g}\left(
			-|\xi|, -\dfrac{\xi}{|\xi|}
			\right)			
	\right),\\
	&\hat{g}(\tau,\theta) = \dfrac{1}{\sqrt{2\pi}}\int\limits_{\R}e^{i\tau s}
	g(s,\theta)ds, \, \tau \in \R, \, \theta\in \Sp^{n-1},
\end{align}
where $g\in C(\Sp^{n-1}, L^2(\R))$, $\xi\in\R^n$. The validity of formulas \eqref{fourier_dual_radon}, \eqref{fourier_dual_derivative} (in the sense of distributions) follows from  Theorem 1.4 of  \cite{natterer1986mathematics}.

\subsection{\large The case of odd $n$}
Using identity \eqref{radon_symmetry} we get
\begin{equation} \label{pre_odd_final_property}
	g(s,\theta) = g(-s,-\theta), \text{ for all } s\in \R, \, \theta\in \Sp^{n-1}, 
\end{equation}
where
\begin{equation}\label{pre_odd_final}
	g(s,\theta) = \left[R_{W_s}f(s,\theta) - R_{w_0}f(s,\theta)\right].
\end{equation}
From \eqref{pre_odd_final_property}, we obtain 
the same symmetry for the Fourier transform $\hat{g}(\cdot, \theta)$ of $g(\cdot,\theta)$:
\begin{align}\nonumber \label{symmetry_odd_fourier_transform}
	\hat{g}(t,\theta) &= \dfrac{1}{\sqrt{2\pi}}\int\limits_{\R}g(s,\theta)e^{its}ds \\
	&=\dfrac{1}{\sqrt{2\pi}}\int\limits_{\R}g(-s,-\theta)e^{i(-s)(-t)}ds\\
	&=\dfrac{1}{\sqrt{2\pi}}\int\limits_{\R}g(s,-\theta)e^{-its}ds = \hat{g}(-t,-\theta),\nonumber
	t\in \R, \, \theta\in \Sp^{n-1}.
\end{align} 
For odd $n$, from identities \eqref{necessity_odd_sym}, \eqref{fourier_dual_derivative} it follows that 
\begin{equation}\label{summ_zero_odd}
	\hat{g}\left(|p|, \dfrac{p}{|p|}\right) + 
	\hat{g}\left(-|p|, -\dfrac{p}{|p|}\right) = 0 \text{ in } L^2_{loc}(\R^n).
\end{equation}
Using \eqref{symmetry_odd_fourier_transform}, \eqref{summ_zero_odd} we obtain
\begin{align}
	&\begin{cases}\label{odd_final}
		\hat{g}\left(
			|p|, \dfrac{p}{|p|}
		\right) = 0, \\
		\hat{g}\left(-|p|, -\dfrac{p}{|p|}\right) = 0
	\end{cases} \Leftrightarrow \hat{g} = 0 \Leftrightarrow g = 0.
\end{align}
Formula \eqref{radon_eq} for odd $n$ follows from \eqref{pre_odd_final}, \eqref{odd_final}.
\subsection{The case of even $n$}
We consider
\begin{equation}\label{even_formula_apply}
	g(s,\theta) = \Hb\left[
		R_{W_s}f - R_{w_0}f
	\right](s,\theta), \, s\in \R, \, \theta\in \Sp^{n-1},
\end{equation}
arising in \eqref{necessity_even_sym}. Using the identity
\begin{align}
	& \mathbb{H}\left[
		R_{W_s}f - R_{w_0}f
	\right](-s,-\theta) = \dfrac{1}{\pi} p.v. \int\limits_{\R}
	\dfrac{R_{W_s}f(t,-\theta) - R_{w_0}f(t,-\theta)}{-s-t}dt\\
	&= \dfrac{1}{\pi} p.v. \int\limits_{\R} \nonumber
	\dfrac{R_{W_s}f(-t,-\theta) - R_{w_0}f(-t,-\theta)}{-s +t}dt \\
	& = -\dfrac{1}{\pi}p.v.\int\limits_{\R} \nonumber
	\dfrac{R_{W_s}f(t,\theta) - R_{w_0}f(t,\theta)}{s-t}dt = -
	\mathbb{H}\left[
		R_{W_s}f - R_{w_0}
	\right](s,\theta),\,
\end{align}
we obtain
\begin{equation}\label{even_g_symmetry}
	g(s,\theta) = -g(-s,-\theta), \text{ for all }s\in \R, \, \theta\in \Sp^{n-1}.
\end{equation}
From \eqref{even_g_symmetry}, similarly with \eqref{symmetry_odd_fourier_transform}, we obtain the 
same symmetry for the Fourier transform $\hat{g}(\cdot,\theta)$ of $g(\cdot,\theta)$:
\begin{align} \label{symmetry_even_fourier_transform}
	\hat{g}(t,\theta) = -\hat{g}(-t,-\theta), \, t\in \R, \theta\in \Sp^{n-1}.
\end{align} 
For $n$ even, from the property of the Hilbert transform
\begin{align}
	\Hb\left[\phi^{(k)}\right] = \left(\Hb \left[\phi\right]\right)^{(k)}, \, 
	\phi\in C_0^k(\R), \nonumber
\end{align}
where this identity holds in the sense of distributions if $\phi\in C_0(\R)$,
and identities 
\eqref{necessity_even_sym}, \eqref{fourier_dual_derivative} it follows that
\begin{equation}\label{odd_summ_zero}
	\hat{g}\left(|p|, \dfrac{p}{|p|}\right) -
	\hat{g}\left(-|p|, -\dfrac{p}{|p|}\right) = 0 \text{ in } L^2_{loc}(\R^n).
\end{equation}
Using \eqref{symmetry_even_fourier_transform}, \eqref{odd_summ_zero} we again obtain \eqref{odd_final} but already for even $n$.
Due to \eqref{odd_final}, \eqref{even_formula_apply} we have 
\begin{equation}
	\Hb\left[\label{pre_pre_even_final}
		R_{W_s}f - R_{w_0}f
	\right] = 0.
\end{equation}
Formula \eqref{radon_eq} for even $n$ follows from \eqref{pre_pre_even_final}, invertibility of the Hilbert transform on $L^p, \, p > 1$ and the fact that $R_Wf\in C_0(\R\times \Sp^{n-1})$.
\par Lemma 1 is proved.

\renewcommand{\thesection}{\large 6}
\section{\large Proof of Lemma 2}
Suppose that 
\begin{equation}
	W_s(y,\theta) - w_0(y) = z \neq 0
\end{equation}
for some $y\in \R^n, \, \theta\in \Sp^{n-1}, \, z\in \mathbb{C}$. Since $W$ satisfies \eqref{weight_class}, then for any $\varepsilon > 0$ there exists $\delta(\varepsilon) > 0$ such that 
\begin{equation}\label{weight_continuity}
	\forall \, y' : |y' - y| < \delta \rightarrow
	|W_s(y',\theta) - w_0(y') - z| < \varepsilon,
\end{equation}
for fixed $y,\theta$. \\
\indent Let $f\in C_0(\R^n), \, f\geq 0$ and satisfies the conditions
\begin{align}
	&f(y')\equiv 1, \,  y' \in B_{\delta/2}(y),\\
	&supp\,f \subset B_{\delta}(y),
\end{align}
where $B_\delta(y)$ is the open ball with radius $\delta$, centered at $y$, $\delta = \delta(\varepsilon), \, 0 < \varepsilon < |z|$. It suffices to show that
\begin{equation}\label{lemma_2_contradiction}
	|R_{W_s}f(y\theta,\theta) - R_{w_0}f(y\theta, \theta) | > 0,
\end{equation}
which contradicts the condition of the lemma.
\par The identity \eqref{lemma_2_contradiction} follows from the formulas
\begin{align}\nonumber
	&|R_{W_s}(y\theta,\theta) - R_{w_0}(y\theta,\theta)| = \Bigl|\int\limits_{x\theta = y\theta}f(x)(W_s(x,\theta) - w_0(x))dx\Bigr|\\
	&=\Bigl|\nonumber
		\int\limits_{x\theta=y\theta}f(x)(W_s(x,\theta)-w_0(x)-z)dx + z\int\limits_{x\theta=y\theta}f(x)dx
	\Bigr|\\ \nonumber
	&\geq |z|\int\limits_{x\theta = y\theta}f(x)dx - 
		\int\limits_{x\theta = y\theta}f(x)\left|W_s(x,\theta) - w_0(x)-z\right|dx \\
	&\geq (|z| - \varepsilon) \int\limits_{x\theta=y\theta}f(x)dx > 0, \text{ for } 0< \varepsilon < |z|.
\end{align}
\indent Lemma 2 is proved.

%
%

\begin{thebibliography}{99}
\thispagestyle{myheadings}

\bibitem{boman1987support}
J.~Boman and E.~T. Quinto,
\newblock {\em Support theorems for real-analytic {R}adon transforms},
\newblock Duke Mathematical J., 55(4), (1987), 943-948.

\bibitem{boman2004novikov}
J.~Boman and J.~Str{\"o}mberg,
\newblock Novikov's inversion formula for the attenuated {R}adon transform -- a new approach,
\newblock {\em The Journal of Geometric Analysis}, 14(2), (2004), 185--198.

\bibitem{chang1978method}
L.~Chang,
\newblock A method for attenuation correction in radionuclide computed
  tomography,
\newblock {\em IEEE Transactions on Nuclear Science}, 25(1), (1978), 638--643.

\bibitem{deans2007radon}
S.~R. Deans,
\newblock The {R}adon {T}ransform and some of {I}ts Applications,
\newblock {\em Courier Corporation}, 2007.

\bibitem{durrani1984radon}
T.~Durrani and D.~Bisset,
\newblock The {R}adon transform and its properties,
\newblock {\em Geophysics}, 49(8), (1984), 1180--1187.

\bibitem{gel2014integral}
I.~M. Gel'fand, M.~I. Graev, and N.~Ya. Vilenkin,
\newblock Integral {G}eometry and {R}epresentation {T}heory,
\newblock {\em Academic press,} Vol. 5, 2014.

\bibitem{gindikin2010remark}
S.~Gindikin,
\newblock A remark on the weighted {R}adon transform on the plane,
\newblock {\em Inverse Problems and Imaging}, 4, (2010), 649--653.

\bibitem{grangeat1991mathematical}
P.~Grangeat,
\newblock Mathematical framework of cone beam 3{D} reconstruction via the first
  derivative of the {R}adon transform,
\newblock {\em Mathematical methods in tomography, Springer}, (1991), 66-97.

\bibitem{helgason2011radon}
S.~Helgason,
\newblock The Radon Transform on {$R^n$},
\newblock {\em Springer}, 2011.

\bibitem{kunyansky1992generalized}
L.~A. Kunyansky,
\newblock Generalized and attenuated {R}adon transforms: restorative approach
  to the numerical inversion,
\newblock {\em Inverse Problems}, 8(5), (1992), 809-819.

\bibitem{ludwig1966radon}
D.~Ludwig,
\newblock The {R}adon transform on {E}uclidean space,
\newblock {\em Communications on Pure and Applied Mathematics}, 19(1), (1996), 49--81.

\bibitem{natterer1986mathematics}
F.~Natterer,
\newblock The {M}athematics of {C}omputerized {T}omography, Vol. 32,
\newblock {\em SIAM}, 1986.

\bibitem{novikov2002inversion}
R.~G. Novikov,
\newblock An inversion formula for the attenuated {X}-ray transformation,
\newblock {\em Arkiv f{\"o}r matematik}, 40(1), (2002), 145--167.

\bibitem{novikov2011weighted}
R.~G. Novikov,
\newblock Weighted {R}adon transforms for which {C}hang's approximate inversion
  formula is exact,
\newblock {\em Russian Mathematical Surveys}, 66(2), (2011), 442--443.

\bibitem{radon1994uber}
{J}. {R}adon,
\newblock {U}ber die {B}estimmung von {F}unktionen durch ihre {I}ntegralwerte
  l{\"a}ngs gewisser {M}annigfaltigkeiten [{O}n the determination of functions by
  their integral values along certain manifolds], 75 years of {R}adon transform
  ({V}ienna, 1992), conf.
\newblock {\em Proc. Lecture Notes Math. Phys., IV, Int. Press, Cambridge, MA},
   (1994), 324--339.

\bibitem{tretiak1980exponential}
O.~Tretiak and C.~Metz,
\newblock The exponential {R}adon transform,
\newblock {\em SIAM Journal on Applied Mathematics}, 39(2), (1980), 341--354.
\end{thebibliography}

\end{document}